\documentclass[11pt,a4paper,fullpage]{article}

\usepackage{amsmath,amsfonts,amssymb,amsthm,graphics,amscd}

\newtheorem{Teorema}{Theorem}[section]

\newtheorem{Posledica}[Teorema]{Corollary}
\newtheorem{Primer}[Teorema]{Example}
\newtheorem{Lema}[Teorema]{Lemma}
\newtheorem{Primedba}[Teorema]{Remark}

\newcommand\scalemath[2]{\scalebox{#1}{\mbox{\ensuremath{\displaystyle #2}}}} 
\numberwithin{equation}{section}

\def\N{\mathcal{N}}
\def\R{\mathcal{R}}
\def\B{\mathcal{B}}

\begin{document}
	\title {Completion problem of upper triangular $3\times3$ operator matrices on arbitrary Banach spaces}

\author{Nikola Sarajlija\footnote{corresponding author: Nikola Sarajlija; University of Novi Sad, Faculty of Sciences, Novi Sad 21000, Serbia; {\it e-mail}: {\tt nikola.sarajlija@dmi.uns.ac.rs}} \ and\ Dragan S. Djordjevi\'c\footnote{University of Ni\v s, Faculty of Sciences and Mathematics, Ni\v s 18000, Serbia; {\it e-mail}: {\tt dragandjordjevic70@gmail.com}} \ \footnote{Authors are supported by the Ministry of Education, Science and Technological Development of the Republic of Serbia under grants no. 451-03-68/2022-14/200124, 451-03-137/2025-03/ 200125 and 451-03-136/2025-03/ 200125} }
\maketitle

\begin{abstract}
We solve a completion problem of a $3\times3$ upper triangular operator matrix acting on a direct sum of Banach spaces and hence generalize the famous result of J. K. Han, H. Y. Lee and W. Y. Lee (Proc. Amer. Math. Soc. \textbf{128} (1) (2000), 119-123) to a greater dimension of a matrix. Our main tools are Harte's ghost of an index theorem and  Banach spaces embeddings. We overcome the lack of orthogona\-li\-ty in Banach spaces by exploiting decomposition properties of inner regular operators. We provide some illustrative examples, and we comment on further generalizations to matrix dimensions $n>3$.
\end{abstract}

\textit{$2020$ Math. Subj. Class:} 47A08, 47A05, 47A53.

\vspace{2mm}

\section{Introduction and preliminaries}

Throughout this text, let $X,Y,Z$ be arbitrary  Banach spaces. Let $\mathcal{B}(X,Y)$ stand for the collection of all bounded linear operators from $X$ to $Y$. By agreement, $\mathcal{B}(X)=\mathcal{B}(X,X)$. If $T\in\mathcal{B}(X,Y)$, then $\mathcal{N}(T)$ and $\mathcal{R}(T)$ stand for the null and the range space of $T$, respectively. 

Let $T\in\mathcal{B}(X,Y).$ We put $\alpha(T)=\dim\mathcal{N}(T)$ and $\beta(T)=\dim Y/\mathcal{R}(T)$. In the previous sentence, $\dim$ is the algebraic dimension of a subspace. If $X$ is a Hilbert space, we use the notation $\dim_hX$ for an orthogonal dimension of a Hilbert space $X$ (term 'orthogonal dimension' refers to the cardinality of any orthonormal bases for an underlying Hilbert space). Different types of regularities can be described using notions $\alpha(\cdot)$ and $\beta(\cdot)$. For example:\\
-$T$ is a left Fredholm operator if $\alpha(T)<\infty$ and $\mathcal{R}(T)$ is closed and complemented in $Y$;\\
-$T$ is a right Fredholm operator if $\beta(T)<\infty$ and $\mathcal{N}(T)$ is complemented in $X$;\\
-$T$ is a Fredholm operator if $\alpha(T),\beta(T)<\infty$;\\
-$T$ is left invertible if $T$ is left Fredholm with $\alpha(T)=0$;\\
-$T$ is right invertible if $T$ is right Fredholm  with $\beta(T)=0$;\\
-$T$ is invertible if $T$ is Fredholm invertible with $\alpha(T)=\beta(T)=0$.\\
Notice that range closedness is not required in the definition of a Fredholm operator. However, this is implicitly included in that definition, since a well known result states that $\beta(T)<\infty$ always implies range closedness of $T$ (see for example Lemma 4.38 in \cite{ABRAMOVICH}).

This article deals with invertibility properties of  $3\times3$ upper triangular block operators whose diagonal entries are known. We denote such an operator by $M_{D,E,F}$, that is
\begin{equation}\label{OSNOVNA}
M_{D,E,F}=\begin{pmatrix}A & D & E\\ 0 & B & F\\ 0 & 0 & C\end{pmatrix},
\end{equation}
where $D,E,F$ are variables. It is understood that $M_{D,E,F}$ acts on the space $X\oplus Y\oplus Z$, while its entries act on appropriate domains: $A\in\mathcal{B}(X)$, $D\in\mathcal{B}(Y,X)$ and so on.

Upper triangular operator matrices of dimension 2 are well studied in terms of their invertibility properties (see \cite{OPERATORTHEORY}, \cite{DU}, \cite{HAN}, \cite{KOLUNDZIJA} etc.). However, research on invertibility of upper triangular operators of dimension 3 has been neglected until a few years ago (\cite{DONG},\cite{WU3}). We show that the well-known result of J. K. Han, H. Y. Lee and W. Y. Lee \cite[Theorem 2]{HAN} can be extended to the case of $3\times3$ upper triangular operator matrices. Let us also mention that there has been some interest in the investigation of $n\times n$ upper triangular operator matrices for general $n>3$, but we shall not pursue this point here (see for example \cite{2016},\cite{SARAJLIJA2} and the last section of this article). 

In \cite{HAN} the authors exploited decomposition properties of inner regular operators (see (0.2) in \cite{HAN}), and we pursue such an idea. The class of inner regular operators consists of operators $T\in\mathcal{B}(X,Y)$ that can be expressed in the form $T=TT'T$ for some $T'\in\mathcal{B}(Y,X)$.  It is well-known that $T\in\mathcal{B}(X,Y)$ is inner regular if and only if its kernel and range are closed and complemented subspaces \cite[Corollary 1.1.5]{DJORDJEVIC}. For example, every (left, right) Fredholm operator is inner regular, and every (left, right) invertible operator is inner regular. Instead of inner regular we will usually say regular for short.

We will make use of the following definition introduced in \cite{OPERATORTHEORY}: we say that a Banach space $X$ can be embedded in a Banach space $Y$, denoted by $X\preceq Y$, provided that there exists a left invertible operator $A\in\B(X,Y)$. Then, it is obvious that $X\cong Y$ if and only if $X\preceq Y$ and $Y\preceq X$. If $X,Y$ are Hilbert spaces, then $X\preceq Y$ if and only if $\dim_hX\leq\dim_hY$. If $U$ is a closed subspace of a Banach space $V$,  we will sometimes use the following notation for the quotient space: $\dfrac{V}{U}=V/U$.

In the sequel, we will find a huge benefit of the following matrix decomposition:
\begin{equation}\label{RAZLAGANJE}
M_{D,E,F}=\scalemath{0.9}{\begin{pmatrix}I & 0 & 0\\0 & I & 0\\0 & 0 & C\end{pmatrix}\begin{pmatrix}I & 0 & E\\0 & I & F\\0 & 0 & I\end{pmatrix}\begin{pmatrix}I & 0 & 0\\0 & B & 0\\0 & 0 & I\end{pmatrix}\begin{pmatrix}I &D & 0\\0 & I & 0\\0 & 0 & I\end{pmatrix}\begin{pmatrix}A & 0 & 0\\0 & I & 0\\0 & 0 & I\end{pmatrix}}
\end{equation}
This representation has already been used in \cite{DONG}. Notice that the second and the fourth factor in (\ref{RAZLAGANJE}) are invertible matrices for all $D\in\mathcal{B}(Y,X),E\in\mathcal{B}(Z,X),F\in\mathcal{B}(Z,Y)$.\\

The following results will be used  in the proof of our main theorem.
\begin{Lema}\label{JEDNAKOSTI}
Let $S,T\in\mathcal{B}(X)$. If $T$ is invertible, then:
\begin{enumerate}
\item[{\rm1)}] $\mathcal{R}(TS)\cong\mathcal{R}(S)$;
\item[{\rm2)}] $\mathcal{R}(ST)=\mathcal{R}(S)$;
\item[{\rm3)}] $\mathcal{N}(ST)\cong\mathcal{N}(S)$;
\item[{\rm4)}]$\mathcal{N}(TS)=\mathcal{N}(S)$.
\end{enumerate}
\end{Lema}
\begin{proof} Parts $2)$ and $4)$ are straightforward. In $1)$ and $3)$ the desired isomorphism is $T$.
\end{proof}

\begin{Lema}\label{JACAINVERTIBILNOST}
Consider $M_{D,E,F}$ and its diagonal operators $A,B,C$. If any three of those four operators are invertible for every $D\in\mathcal{B}(Y,X),E\in\mathcal{B}(Z,X),F\in\mathcal{B}(Z,Y)$, then the fourth is invertible as well.
\end{Lema}
\begin{proof}This is obvious from (\ref{RAZLAGANJE}). \end{proof}

The following lemma is well-known in the literature (see for example \cite[Lemma 2.3]{OPERATORTHEORY}).
\begin{Lema}\label{LEMA2.3}
If $X,Y,Z$ are Banach spaces then
$$(X\times Y\cong X\times Z\ \text{and} \ \dim X<\infty) \Rightarrow\ Y\cong Z.$$
\end{Lema}



The following result is known as Harte's ghost of an index theorem (see \cite{HARTE} and \cite{HARTE2}).
\begin{Teorema}\label{ghost}
If $T\in\B(X,Y)$, $S\in \B(Y,Z)$ and $ST\in\B(X,Z)$ are regular, then
 $$\N(T)\times\N(S)\times Z/\R(ST)\cong \N(ST)\times Y/\R(T)\times Z/\R(S).$$
\end{Teorema}

\section{Characterization of invertibility of $M_{D,E,F}$}

In this section we provide necessary and sufficient conditions for invertibility of $M_{D,E,F}$ (Theorem \ref{GLAVNA}), and afterwards consider under what a priori assumption these conditions become equivalent (Corollary \ref{GLAVNA2}, Theorem \ref{hil}). 

We prove the following theorem.

\begin{Teorema}\label{GLAVNA} Let $X,Y,Z$ be Banach spaces, and let $A\in\B(X)$, $B\in\B(Y)$, $C\in\B(Z)$. Assume that $B$ is regular.  Consider the following statements:
\begin{enumerate}
\item[{\rm1)}]
\begin{enumerate}
\item[{\rm a)}] $A$ is left invertible and $C$ is right invertible;
\item[{\rm b)}] $\mathcal{N}(B)\preceq X/\mathcal{R}(A)$ and $Y/\mathcal{R}(B)\preceq\mathcal{N}(C)$;
\item[{\rm c)}] $\dfrac{X/\mathcal{R}(A)}{\mathcal{R}(J_1)}\cong\dfrac{\mathcal{N}(C)}{\mathcal{R}(J_2)}$ for left invertible operators $J_1:\N(B)\to X/\R(A)$ and $J_2:Y/\R(B)\to\N(C)$ which realize relations $\preceq$ in {\rm 1) b)}.
\end{enumerate}
\item[{\rm 2)}] There exist $D\in\mathcal{B}(Y,X),E\in\mathcal{B}(Z,X),F\in\mathcal{B}(Z,Y)$ such that $M_{D,E,F}$ is invertible.
\item[{\rm 3)}] 
\begin{enumerate}
\item[{\rm a)}] $A$ is left invertible and $C$ is right invertible;
\item[{\rm b)}] $\mathcal{N}(B)\preceq X/\mathcal{R}(A)$ and $Y/\mathcal{R}(B)\preceq\mathcal{N}(C)$;
\item[{\rm c)}] $\mathcal{N}(B)\times\mathcal{N}(C)\cong X/\mathcal{R}(A)\times Y/\mathcal{R}(B).$
\end{enumerate}
\end{enumerate}
Then $1)\Rightarrow 2) \Rightarrow 3)$.
\end{Teorema}

\begin{proof} 
1)$\implies$2): Suppose that 1) holds. Since $A,B$ and $C$ are inner regular operators, there exist closed subspaces: $X_1$ of $X$, $Y_1$ and $Y_2$ of $Y$, and $Z_1$ of $Z$ such that:
$$X=X_1\oplus \R(A),\quad Y=Y_1\oplus \R(B),\quad Y=Y_2\oplus\N(B),\quad Z=Z_1\oplus \N(C).$$
Consequently, 
$$X/\R(A)\cong X_1,\quad  Y/\R(B)\cong Y_1, \quad Y/\N(B)\cong Y_2,\quad Z/\N(C)\cong Z_1.$$

Condition 1) b) implies the existence of left invertible mappings $J_1:\mathcal{N}(B)\rightarrow X_1$ and $J_2:Y_1\rightarrow\mathcal{N}(C)$.  Consider their invertible reductions  $J_1:\N(B)\to \R(J_1)$ and $J_2:Y_1\to \R(J_2)$, which are denoted by the same symbols.
There exist closed subspaces $\R(J_1)'$ and $\R(J_2)'$ such that 
$$X_1=\R(J_1)'\oplus\R(J_1),\quad \N(C)=\R(J_2)'\oplus\R(J_2).$$ Since $X_1\cong X/\R(A)$, notice that we are now able to summon condition 1) c). Namely, this condition implies that there exists an isomorphism $J:\R(J_2)'\to\R(J_1)'$. For notational consistency we denote $M_1=\mathcal{R}(J_1),M_1'=\mathcal{R}(J_1)'$ and $M_2=\mathcal{R}(J_2),M_2'=\mathcal{R}(J_2)'$.

Define 
$$D=\begin{pmatrix}0 &0\\ 0 & 0\\ 0 & J_1\end{pmatrix}:\begin{pmatrix}Y_2\\ \mathcal{N}(B)\end{pmatrix}=Y\rightarrow X=\begin{pmatrix}\mathcal{R}(A)\\ M_1'\\ M_1\end{pmatrix},$$  
$$F=\begin{pmatrix}0 &0 & 0\\0 & 0 &  J_2^{-1}\end{pmatrix}:\begin{pmatrix}Z_1\\ M_2'\\ M_2\end{pmatrix}=Z\rightarrow Y=\begin{pmatrix}\mathcal{R}(B)\\ Y_1\end{pmatrix},$$ 
and
$$E=\begin{pmatrix}0 & 0 & 0\\ 0 & J &0\\ 0 & 0 & 0\end{pmatrix}:\begin{pmatrix}Z_1\\ M_2'\\ M_2\end{pmatrix}=Z\rightarrow X=\begin{pmatrix}\mathcal{R}(A)\\ M_1'\\ M_1\end{pmatrix}.$$ 
  Since $J_1$, $J_2^{-1}$ and $J$ are isomorphisms between appropriate subspaces, it is obvious that $D\in\mathcal{B}(Y,X),E\in\mathcal{B}(Z,X),F\in\mathcal{B}(Z,Y)$.
  
  To prove that $M_{D,E,F}$ is invertible, notice the following equalities. 
  We have
  $$\aligned &
  \R(M_{D,E,F})=\R\left(
  \bmatrix A\\0\\0\endbmatrix\right)+
  \R\left(\bmatrix D\\B\\0\endbmatrix\right)+\R\left(\bmatrix E\\F\\C\endbmatrix\right)\\
  &=\big(\R(A)+M_1+M_1'\big)+\big(\R(B)+Y_1\big)+\R(C)=X\oplus Y\oplus Z,\endaligned
  $$
so $M_{D,E,F}$ is onto. 

Moreover, if $w=\bmatrix x\\y\\z\endbmatrix\in X\oplus Y\oplus Z$ and $M_{D,E,F}w=0$, we have
$$Ax+Dy+Ez=0,\quad By+Fz=0,\quad Cz=0.$$
From $Cz=0$ we get $z\in\N(C)=M_2'\oplus M_2$. We know that $By\in\R(B)$ and $Fz\in Y_1$. Thus, from $By+Fz=0$ we get $By=0$ and $Fz=0$. Hence, $y\in\N(B)$ and $z\in M_2'$. We have $Ax\in\R(A)$, $Dy\in D(\N(B))=J_1(\N(B))=\R(J_1)=M_1$ and $Ez\in E(\R(J_2)')=J(\R(J_2)')=\R(J_1)'=M_1'$. Hence, from $Ax+Dy+Ez=0$ we conclude $Ax=0$, $Dy=J_1y=0$ and $Jz=0$, implying that $x=0$, $y=0$ and $z=0$. Thus, $M_{D,E,F}$ is one-to-one.

\vspace{3mm}

2)$\implies$3): Assume that $M_{D,E,F}$ is invertible for some $D$, $E$ and $F$ defined on appropriate domains. Consider factorization (\ref{RAZLAGANJE}) to conclude that $A$ is left invertible and $C$ is right invertible, thus the condition 3) a) follows.

Denote the product of the first two factors in (\ref{RAZLAGANJE}) by $S$, the product of the last three factors by $T$, and apply Theorem \ref{ghost}. We get
\begin{equation*}
\mathcal{N}(S)\times\mathcal{N}(T)\times\lbrace0\rbrace\cong\lbrace0\rbrace\times\begin{pmatrix}X\\ Y\\ Z\end{pmatrix}/\mathcal{R}(S)\times\begin{pmatrix}X\\ Y\\ Z\end{pmatrix}/\mathcal{R}(T),
\end{equation*}
where we have used invertibility of $ST=M_{D,E,F}$.

Next, in view of Lemma \ref{JEDNAKOSTI}, the previous congruence becomes
\begin{equation*}
\mathcal{N}(C)\times\mathcal{N}(T)\cong Z/\mathcal{R}(C)\times\begin{pmatrix}X\\ Y\\ Z\end{pmatrix}/\mathcal{R}(T). 
\end{equation*}
Observe that $T$ is left invertible, hence $\mathcal{N}(T)=\lbrace0\rbrace$, and that $C$ is right invertible, hence $Z/\mathcal{R}(C)=\lbrace0\rbrace$. Therefore, we finally get
\begin{equation}\label{PRVA}
\mathcal{N}(C)\cong\begin{pmatrix}X\\ Y\\ Z\end{pmatrix}/\mathcal{R}(T). 
\end{equation}

If we denote the product of the first three factors in (\ref{RAZLAGANJE}) by $S'$, and the product of the last two by $T'$, in the similar manner one can get
\begin{equation}\label{STOSEDMA}
\mathcal{N}(S')\cong X/\mathcal{R}(A).
\end{equation}
Relations (\ref{PRVA}) and (\ref{STOSEDMA}) immediately imply condition 3) b). Namely, $\mathcal{N}(B)\cong\mathcal{N}\left(\begin{pmatrix}I & 0 & 0\\ 0 & B & 0\\ 0 & 0 & I\end{pmatrix}\right)\subseteq\mathcal{N}(S')\cong X/\mathcal{R}(A)$ by definition of $S'$ and \eqref{STOSEDMA}, and similar for $Y/\mathcal{R}(B)\preceq\mathcal{N}(C)$.

Now, denote the product of the second and the third factor in (\ref{RAZLAGANJE}) by $\widetilde{S}$, the product of the last two factors by $\widetilde{T}$, and again apply Theorem \ref{ghost}  to $\widetilde{S}$ and $\widetilde{T}$. We get
\begin{equation*}
\mathcal{N}(\widetilde{S})\times\mathcal{N}(\widetilde{T})\times\begin{pmatrix}X\\ Y\\ Z\end{pmatrix}/\mathcal{R}(\widetilde{S}\widetilde{T})\cong\mathcal{N}(\widetilde{S}\widetilde{T})\times\begin{pmatrix}X\\ Y\\ Z\end{pmatrix}/\mathcal{R}(\widetilde{S})\times\begin{pmatrix}X\\ Y\\ Z\end{pmatrix}/\mathcal{R}(\widetilde{T}).
\end{equation*}

Next, in view of Lemma \ref{JEDNAKOSTI}, the previous congruence becomes
\begin{equation*}
\mathcal{N}(B)\times\mathcal{N}(A)\times\begin{pmatrix}X\\ Y\\ Z\end{pmatrix}/\mathcal{R}(\widetilde{S}\widetilde{T})\cong\mathcal{N}(\widetilde{S}\widetilde{T})\times Y/\mathcal{R}(B)\times X/\mathcal{R}(A).
\end{equation*}
Observe that $A$ and $\widetilde{S}\widetilde{T}$ are left invertible, hence their null spaces are zero. Therefore, we finally get
\begin{equation}\label{DRUGA}
\mathcal{N}(B)\times\begin{pmatrix}X\\ Y\\ Z\end{pmatrix}/\mathcal{R}(\widetilde{S}\widetilde{T})\cong Y/\mathcal{R}(B)\times X/\mathcal{R}(A).
\end{equation}

Now, there is one more use of Lemma \ref{JEDNAKOSTI}  to conclude that (\ref{PRVA}) and (\ref{DRUGA}) imply $\mathcal{N}(B)\times\mathcal{N}(C)\cong X/\mathcal{R}(A)\times Y/\mathcal{R}(B),$ which is 3) c). 
\end{proof}


Naturally, one would like to know if there are regular operators $B$ such that conditions 1) and 3) in Theorem \ref{GLAVNA} become equivalent. Implication $1)\Rightarrow3)$ always holds, and with regards to Lemma \ref{LEMA2.3} it is obvious that opposite is true for the class of Fredholm operators. The following theorem is our main result. It is at the same time extension of \cite[Theorem 2]{HAN} to matrix dimension $n=3$, generalization of \cite[Corollary 2.4]{WU3} to the Banach space setting, and improvement of \cite[Theorem 2.13]{SARAJLIJA3} for $n=3$.

\begin{Posledica}\label{GLAVNA2}
Assume that $B\in\mathcal{B}(Y)$ is Fredholm. Then the following statements are equivalent: 
\begin{enumerate}

\item[{\rm1)}] There exist $D\in\mathcal{B}(Y,X),E\in\mathcal{B}(Z,X),F\in\mathcal{B}(Z,Y)$ such that $M_{D,E,F}$ is invertible.

\item[{\rm2)}] \begin{enumerate}
\item[{\rm a)}] $A$ is left invertible and $C$ is right invertible;
\item[{\rm b)}]  $\mathcal{N}(B)\preceq X/\mathcal{R}(A)$ and $Y/\mathcal{R}(B)\preceq\mathcal{N}(C)$;
\item[{\rm c)}]  $\mathcal{N}(B)\times\mathcal{N}(C)\cong X/\mathcal{R}(A)\times Y/\mathcal{R}(B)$.
\end{enumerate}
\end{enumerate}
\end{Posledica}

Previous corollary is  very reminiscent to the statement of \cite[Theorem 2]{HAN}.

\begin{Primer}
    One might consider conditions stated in Theorem \ref{GLAVNA} and Corollary \ref{GLAVNA2} cumbersome and hard-to-work with. However, this is not the case. For example, take $X=Y=Z=l^2$ and let $B\in\mathcal{B}(Y)$ be any diagonal operator on $l^2$ of the form $T(x_1,x_2,...)=(\lambda_1x_1,\lambda_2x_2,...)$ with $sup\vert\lambda_i\vert<\infty$ and $min\vert\lambda_i\vert>0$. Then operator $B$ is Fredholm. Now, among all left invertible operators on $\mathcal{B}(X)$ choose $A$ defined as $A(x_1,x_2,...)=(0,x_1,0,x_2,...)$. Finally, among all right invertible operators on $\mathcal{B}(Z)$ choose $C(x_1,x_2,x_3,x_4,...)=(x_1+x_2,x_3+x_4,...)$. One can easily check that condition 2) in Corollary 2.2 is fulfilled. Notice that the same example works for the statements to follow, as well.
\end{Primer}

We prove the following result for Hilbert space operators.

\begin{Teorema}\label{hil} Let $X,Y,Z$ be separable Hilbert spaces, $A\in\B(X)$ is left invertible, $B\in\B(Y)$ is inner regular, $C\in\B(Z)$ is right invertible, 
$$\dim_h\N(B)\leq\dim_hX/\R(A)\quad\text{and}\quad \dim_hY/\R(B)\leq\dim_h\N(C).$$
Then the following statements are equivalent;
\begin{enumerate}
    \item[{\rm1)}] $\dfrac{X/\R(A)}{\R(J_1)}\cong\dfrac{\N(C)}{\R(J_2)}$ for some left invertible operators $J_1:\N(B)\to X/\R(A)$ and $J_2:Y/\R(B)\to\N(C)$.
    
    \item[{\rm2)}] $\dim_h\N(B)+\dim_h\N(C)=\dim_hX/\mathcal{R}(A)+\dim_hY/\mathcal{R}(B)$.
\end{enumerate}
\end{Teorema}

\begin{proof}
It is enough to prove implication 2)$\implies$1). Suppose that 2) holds. Left invertible operators  $J_1:\N(B)\to X/\R(A)$ and $J_2:Y/\R(B)\to\N(C)$ exist by the main assumption of this theorem. We have to prove that $J_1$ and $J_2$ can be adjusted  such that 1) is also satisfied. 

We consider several cases and subcases.
\vspace{2mm}

{\it Case I}.  $\dim_h\N(B)<\dim_h X/\R(A)$ and $\dim_h\N(C)\leq\dim_h X/\R(A)$.
\vspace{1mm}

\qquad {\it Subcase I.1.} $X/\R(A)$ is infinite dimensional.

Since 
$$\dim_hY/\R(B)\leq\dim_h\N(C)\leq\dim_hX/\R(A),$$ 
by 2) it follows that $\dim_h\N(C)=\dim_hX/\R(A)$. Then 
$$\dim_hJ_1(\N(B))=\dim_h\N(B)<\dim_h X/\R(A).$$ Thus 
$$\dfrac{X/\R(A)}{\R(J_1)}\cong X/\R(A).$$
Since 
$$\dim_h\R(J_2)=\dim_hY/\R(B)\leq\dim_h X/\R(A)=\dim_h\N(C),$$
we conclude that $J_2$ can be adjusted such that 
$$\dfrac{\N(C)}{\R(J_2)}\cong \N(C)\cong X/\R(A)\cong \dfrac{X/\R(A)}{\R(J_1)}. 
$$
Thus, 1) holds.
\vspace{1mm}

\qquad{\it Subcase I.2.} $X/\R(A)$ is finite dimensional.

Let 
$$k=\dim_h\N(B),\quad l=\dim_h\N(C), \quad  m=\dim_h X/\R(A), \quad n=\dim Y/\R(B).$$ We have
$$k<m,\quad n\leq l\leq m,\quad k+l=m+n, $$
all these quantities are finite, and we get  
$$0<m-k=l-n,$$
which is 1) in finite dimensions.
\vspace{2mm}

{\it Case II}. $\dim_h\N(B)<\dim_h X/\R(A)<\dim_h\N(C)$.
\vspace{1mm}

\qquad{\it Subcase II.1.} $\N(C)$ is infinite dimensional.

We get that
$$\N(B)\times\N(C)\cong\N(C)\text{ and }\dim_h Y/\R(B)=\dim_h\N(C).$$
Since $\dim_h X/\R(A)<\dim_h\N(C)$, for every left invertible $J_1:\N(B)\to X/\R(A)$ is is possible to adjust some left invertible $J_2:Y/\R(B)\to\N(C)$ such that 
$$\dfrac{X/\R(A)}{\R(J_1)}\cong\dfrac{\N(C)}{\R(J_2)}$$
holds.
\vspace{1mm}

\qquad{\it Subcase II.2.} $\N(C)$ is finite dimensional.

Keep $k,l,m,n$ the same as in Subcase I.2. We get
$$k<m<l,\quad n\leq l,\quad l+l=m+n,$$
implying that all these quantities are finite and
$$0<m-k=l-n,$$ which is again 1) in finite dimensions.

\vspace{2mm}

{\it Case III}. $\dim_h \N(B)=\dim_hX/\R(A)$ and $\dim_h\N(C)\leq \dim_hX/\R(A)$.
\vspace{1mm}

\qquad{\it Subcase III.1.} $X/\R(A)$ is infinite dimensional.

From 
$$\aligned
\dim_h J_2(Y/\R(B))&=\dim_h Y/\R(B)\leq\dim_h\N(C)\leq\dim_h X/\R(A)\\
&=\dim_h\N(B)\endaligned$$
we get that for every left invertible $J_2:Y/\R(B)\to\N(C)$ we can find a left invertible $J_1:\N(B)\to X/\R(A)$ such that

$$\dfrac{X/\R(A)}{\R(J_1)}\cong \dfrac{\N(C)}{\R(J_2)}.$$
\vspace{1mm}

\qquad{\it Subcase III.2.} $X/\R(A)$ is finite dimensional.

This is proved in the same way as in the previous finite dimensional subcases. 
\vspace{2mm}

{\it Case IV}. $\dim_h\N(B)=\dim_hX/\R(A)<\dim_h\N(C)$.
\vspace{1mm}

\qquad{\it Subcase IV.1.} $\N(C)$ is infinite dimensional.

We get
$$\N(C)\cong \N(B)\times \N(C)\cong X/\R(A)\times Y/\R(B),$$
implying $\N(C)\cong Y/\R(B)$. Thus, for every left invertible $J_1:\N(B)\to X/\R(A)$ we can adjust a left invertible $J_2:Y/\R(B)\to\N(C)$ such that
$$\dfrac{X/\R(A)}{\R(J_1)}\cong\dfrac{\N(C)}{\R(J_2)}.$$
\vspace{1mm}

\qquad{\it Subcase IV.2.} $\N(C)$ is finite dimensional.

Again, this is a routine. 
\end{proof}

We provide a corollary of this final result, and we refer an interested reader to \cite{MOSLEHIAN} for some deeper research on the topic of block operator matrices .

\begin{Posledica}
    Let $X,Y,Z$ be separable Hilbert spaces, $A\in\B(X)$ is left invertible with $codim_h(\mathcal{R}(A))=\infty$, $B\in\B(Y)$ is inner regular, $C\in\B(Z)$ is right invertible with $dim_h(N(C))=\infty$.
Then $$\dfrac{X/\R(A)}{\R(J_1)}\cong\dfrac{\N(C)}{\R(J_2)}$$ for some left invertible operators $J_1:\N(B)\to X/\R(A)$ and $J_2:Y/\R(B)\to\N(C)$.
\end{Posledica}

\begin{Primedba}
    Notice that Banach spaces appearing in Theorem \ref{GLAVNA} and Corollary \ref{GLAVNA2} have linear bases of arbitrary large cardinality. In other words, infinite dimensional Banach spaces that have appeared thus far need not have countable linear bases. Notice that this was not the case in Theorem \ref{hil}, in which we use orthogonal Hilbert bases that are of cardinality $\aleph_0$ or less. We recommend \cite[Section 7]{HALMOS} for deeper understanding.
\end{Primedba}

\section{Extension to higher dimensions}
It is a natural question whether \cite[Theorem 2]{HAN} and Corollary \ref{GLAVNA2} admit a natural generalization from matrix dimensions $n=2$ and $n=3$ to higher matrix dimensions $n>3$. Since Corollary \ref{GLAVNA2} contains an additional apriori assumption (we demand that one of diagonal operators is Fredholm), one must suppose that a mentioned generalization would incorporate more than one appriori assumtions of the same type. Furthermore, regarding the form of Theorem \ref{GLAVNA}, one must conclude that such a generalization would consider appearing more than two left invertible operators $J_i$ similar to those from 1) c) in Theorem \ref{GLAVNA}. Therefore, a logical conclusion is that such a generalization would be very technical and with cumbersome assumptions. This conclusion is further supported by some earlier works of the first author. Namely, statement \cite[Theorem 2.13]{SARAJLIJA3} confirms such a hypothesis.  \\

\textbf{Data availability statement}\\[3mm]
We declare there are no datasets associated with this work.\\

\textbf{Ethics declaration statement}\\[3mm]
Authors declare there are no conflicts of interest associated with this work and that there is no relevant funding other than mentioned which influenced the writing of this manuscript.  \\

\textbf{Acknowledgements}\\[3mm]
We are grateful to anonymous reviewers for their useful advices regarding this paper.

\end{document}